\documentclass[11pt]{amsart}
\usepackage{amssymb,amsmath,amsfonts,latexsym,amsthm,geometry}
\usepackage{amsmath}
\usepackage{amsfonts}
\usepackage{amssymb}
\usepackage{graphicx}
\usepackage{color}
\usepackage[colorlinks=true,linkcolor=blue,citecolor=blue]{hyperref}

\providecommand{\U}[1]{\protect\rule{.1in}{.1in}}
\providecommand{\U}[1]{\protect\rule{.1in}{.1in}}
\providecommand{\U}[1]{\protect\rule{.1in}{.1in}}
\newtheorem{theorem}{Theorem}[section]

\theoremstyle{definition}

\newtheorem{remark}[theorem]{Remark}

\begin{document}
\title[Sharp Hardy--Littlewood inequality for positive multilinear forms]{Sharp anisotropic Hardy--Littlewood inequality for positive multilinear forms}
\author[D. N\'u\~nez]{Daniel N\'{u}\~{n}ez-Alarc\'{o}n}
\address{Departamento de Matem\'{a}ticas\\
\indent Universidad Nacional de Colombia\\
\indent111321 - Bogot\'a, Colombia}
\email{danielnunezal@gmail.com and dnuneza@unal.edu.co}
\author[D. Pellegrino]{Daniel Pellegrino}
\address{Departamento de Matem\'{a}tica \\
Universidade Federal da Para\'{\i}ba \\
58.051-900 - Jo\~{a}o Pessoa, Brazil.}
\email{pellegrino@pq.cnpq.br and dmpellegrino@gmail.com}
\author[D. M. Serrano]{Diana Marcela Serrano-Rodr\'{\i}guez}
\address{Departamento de Matem\'{a}ticas\\
\indent Universidad Nacional de Colombia\\
\indent111321 - Bogot\'{a}, Colombia}
\email{dmserrano0@gmail.com and diserranor@unal.edu.co}

\begin{abstract}
Using elementary techniques, we prove sharp anisotropic Hardy-Littlewood
inequalities for positive multilinear forms. In particular, we recover an
inequality proved by F. Bayart in 2018.

\end{abstract}
\maketitle

\section{Introduction}

\bigskip All along this paper $\mathbb{K}:=\mathbb{R}$ or $\mathbb{C}$,
$X_{p}:=\ell_{p}$ for $1\leq p<\infty$ and $X_{\infty}:=c_{0}$ and, as usual,
we consider $1/\infty:=0$ and $q^{\ast}$ is the conjugate of $q$, i.e.,
$1/q+1/q^{\ast}=1.$ In 1934, Hardy and Littlewood \cite{hardy} proved five
theorems related to summability of bilinear forms. We are interested in the
last one:

\begin{theorem}
\label{999}(See Hardy and Littlewood \cite[Theorem 5]{hardy}) Let
$p,q\in(1,\infty]$ be such that
\[
\frac{1}{p}+\frac{1}{q}<1.
\]
Then%
\[
\left(  \sum_{j_{1}=1}^{\infty}\left(  \sum_{j_{2}=1}^{\infty}A(e_{j_{1}%
},e_{j_{2}})^{q^{\ast}}\right)  ^{\frac{1}{q^{\ast}}\times\frac{1}{1-\left(
\frac{1}{p}+\frac{1}{q}\right)  }}\right)  ^{1-\left(  \frac{1}{p}+\frac{1}%
{q}\right)  }\leq\left\Vert A\right\Vert ,
\]
for all bounded non-negative (i.e., $A(e_{j_{1}},e_{j_{2}})\geq0$ for all
$\left(  j_{1},j_{2}\right)  \in\mathbb{N\times N}$) bilinear forms
$A:X_{p}\times X_{q}\rightarrow\mathbb{K}$.
\end{theorem}

\bigskip The paper of Hardy and Littlewood was revisited in 1981 by
Praciano-Pereira \cite{ppp} and, recently, by several authors (see, for
instance, \cite{araujo} and the references therein). There are still several
subtle open problems regarding the generalization of the Hardy--Littlewood
inequalities to multilinear forms (see, for instance, \cite{laaron}). In 2018,
using a factorization result due to Schep \cite{Schep84}, Bayart
\cite[Proposition 3.1]{bayart} generalized Theorem \ref{999} as follows:

\begin{theorem}
\label{praciano}(See Bayart \cite[Proposition 3.1]{bayart}) Let $m$ be a
positive integer and $p_{1},\dots,p_{m}\in(1,\infty]$ with
\[
\frac{1}{p_{1}}+\cdots+\frac{1}{p_{m}}<1.
\]
Then
\[
\left(  \sum_{j_{1},...,j_{m}=1}^{\infty}A(e_{j_{1}},...,e_{j_{m}})^{\rho
}\right)  ^{1/\rho}\leq\Vert A\Vert
\]
for all bounded non-negative $m$-linear forms $A:X_{p_{1}}\times\cdots\times
X_{p_{m}}\rightarrow\mathbb{K}$ if, and only if,%

\begin{equation}
\rho\geq\frac{1}{1-\left(  \frac{1}{p_{1}}+\cdots+\frac{1}{p_{m}}\right)  }.
\label{bayo}%
\end{equation}

\end{theorem}

\bigskip In the present paper we prove a new generalization of Theorem
\ref{999}, keeping its anisotropic essence. Following the notation introduced
in \cite{laaron}, let us define $\delta^{s_{k},...,s_{m}}$ by the formula%
\[
\delta^{s_{k},...,s_{m}}:=\frac{1}{1-\left(  \frac{1}{s_{k}}+\cdots+\frac
{1}{s_{m}}\right)  },
\]
for all positive integers $m$ and $k=1,...,m.$ When $\frac{1}{s_{k}}%
+\cdots+\frac{1}{s_{m}}\geq1$ it is convenient to define%
\[
\delta^{s_{k},...,s_{m}}:=\infty.
\]
Also, when $q=\infty$, the notation $\left(  \sum\left\vert x_{j}\right\vert
^{q}\right)  ^{1/q}$ shall represent the supremum of $\left\vert
x_{j}\right\vert .$ \bigskip We prove the following:

\begin{theorem}
\label{hl5multi} Let $p_{1},...,p_{m}\in\lbrack1,\infty],$ $q_{1},...,q_{m}%
\in(0,\infty]$ and $m$ be a positive integer. Then, for any bijection
$\sigma:\{1,...,m\}\rightarrow\{1,...,m\}$ we have
\[
\left(  \sum_{j_{\sigma\left(  1\right)  }=1}^{\infty}\left(  \sum
_{j_{\sigma\left(  2\right)  }=1}^{\infty}\cdots\left(  \sum_{j_{\sigma\left(
m\right)  }=1}^{\infty}A(e_{j_{1}},...,e_{j_{m}})^{q_{m}}\right)
^{\frac{q_{m-1}}{q_{m}}}\cdots\right)  ^{\frac{q_{1}}{q_{2}}}\right)
^{\frac{1}{q_{1}}}\leq\left\Vert A\right\Vert
\]
for all bounded non-negative $m$-linear forms $A:X_{p_{1}}\times\cdots\times
X_{p_{m}}\rightarrow\mathbb{K}$ if, and only if
\[
q_{1}\geq\delta^{p_{\sigma\left(  1\right)  },...,p_{\sigma\left(  m\right)
}},q_{2}\geq\delta^{p_{\sigma\left(  2\right)  },...,p_{\sigma\left(
m\right)  }},...,q_{m-1}\geq\delta^{p_{\sigma\left(  m-1\right)  }%
,p_{\sigma\left(  m\right)  }},q_{m}\geq\delta^{p_{\sigma\left(  m\right)  }%
}.
\]

\end{theorem}

\begin{remark}
\bigskip Note that we do not need the hypothesis
\begin{equation}
\frac{1}{p_{1}}+\cdots+\frac{1}{p_{m}}<1. \label{wa}%
\end{equation}
The paper of Hardy and Littlewood and the recent literature just encompasses
the case (\ref{wa}). For bilinear forms, the complementary case, called by
Hardy and Littlewood \ as case of spaces of type $\alpha$ was investigated in
the seminal paper of M. Riesz \cite{riesz}.
\end{remark}

\begin{remark}
Our result recovers Theorem \ref{praciano}. In fact, if $\frac{1}{p_{1}%
}+\cdots+\frac{1}{p_{m}}<1,$ it is clear that $\delta^{p_{1},...,p_{m}}$ is
the biggest exponent and coincides with the optimal exponent given by
(\ref{bayo}); thus the canonical inclusions of $\ell_{p}$ spaces provides the result.
\end{remark}

\section{The proof}

To simplify the notation we will consider $\sigma(j)=j$ for all $j;$ the other
cases are similar.

\textbf{First Case.} $\frac{1}{p_{1}}+\cdots+\frac{1}{p_{m}}<1.$

\bigskip The proof of the direct implication is a consequence of techniques
used in \cite{laaron}. We present the argument for the sake of completeness.
Let us suppose that $p_{1}>1$ (i.e., $\frac{1}{p_{1}}<1$) and that
\[
\left(  \sum_{j_{1}=1}^{\infty}A(e_{j_{1}})^{q_{1}}\right)  ^{\frac{1}{q_{1}}%
}\leq\left\Vert A\right\Vert
\]
for all continuous non-negative linear forms $A:X_{p_{1}}\rightarrow
\mathbb{K}$. For each $n$ consider the continuous non-negative linear form
$A_{n}(x)=\sum\limits_{j=1}^{n}x_{j}.$ By the H\"{o}lder Inequality, we have
\[
\left\Vert A_{n}\right\Vert \leq n^{\frac{1}{p_{1}^{\ast}}}.
\]
On the other hand%
\[
\left(  \sum_{j=1}^{n}A_{n}(e_{j})^{q_{1}}\right)  ^{\frac{1}{q_{1}}}%
=n^{\frac{1}{q_{1}}},
\]
and, since $n$ is arbitrary,
\[
q_{1}\geq p_{1}^{\ast}=\delta^{p_{1}}.
\]
Thus the case $m=1$, is done. Now, let us proceed by induction. Suppose that
the result is valid for $m-1$ and let
\[
\frac{1}{p_{1}}+\cdots+\frac{1}{p_{m}}<1.
\]
Thus%
\[
\frac{1}{p_{2}}+\cdots+\frac{1}{p_{m}}<1
\]
and the induction hypothesis combined with a simple argument of summability
tells us that, if
\[
\left(  \sum_{j_{1}=1}^{\infty}\left(  \sum_{j_{2}=1}^{\infty}\cdots\left(
\sum_{j_{m}=1}^{\infty}A(e_{j_{1}},...,e_{j_{m}})^{q_{m}}\right)
^{\frac{q_{m-1}}{q_{m}}}\cdots\right)  ^{\frac{q_{1}}{q_{2}}}\right)
^{\frac{1}{q_{1}}}\leq\left\Vert A\right\Vert
\]
for all bounded non-negative $m$-linear forms $A:X_{p_{1}}\times\cdots\times
X_{p_{m}}\rightarrow\mathbb{K}$, then
\begin{align*}
q_{2}  &  \geq\delta^{p_{2},...,p_{m}}\\
&  \vdots\\
q_{m-1}  &  \geq\delta^{p_{m-1},p_{m}}\\
q_{m}  &  \geq\delta^{p_{m}}\text{.}%
\end{align*}
So, we must only show that
\[
q_{1}\geq\delta^{p_{1},...,p_{m}}.
\]
For each $n$ consider the continuous non-negative $m-$linear form
$B_{n}:X_{p_{1}}\times\cdots\times X_{p_{m}}\rightarrow\mathbb{K}$ given by%
\[
B_{n}(x^{\left(  1\right)  },...,x^{\left(  m\right)  })=\sum\limits_{j=1}%
^{n}x_{j}^{\left(  1\right)  }x_{j}^{\left(  2\right)  }...x_{j}^{\left(
m\right)  }.
\]
Since
\[
\frac{1}{\delta^{p_{1},\cdots,p_{m}}}+\sum_{k=1}^{m}\frac{1}{p_{k}}=1\text{,}%
\]
we use the H\"{o}lder inequality and obtain%
\begin{align*}
\left\Vert B_{n}\right\Vert  &  =\sup_{\left\Vert x^{\left(  1\right)
}\right\Vert ,\cdots,\left\Vert x^{\left(  m\right)  }\right\Vert \leq
1}\left\vert \sum\limits_{j=1}^{n}x_{j}^{\left(  1\right)  }x_{j}^{\left(
2\right)  }...x_{j}^{\left(  m\right)  }\right\vert \\
&  \leq\sup_{\left\Vert x^{\left(  1\right)  }\right\Vert ,\cdots,\left\Vert
x^{\left(  m\right)  }\right\Vert \leq1}\left(  \prod_{k=1}^{m}\left(
\sum\limits_{j=1}^{n}\left\vert x_{j}^{\left(  k\right)  }\right\vert ^{p_{k}%
}\right)  ^{1/p_{k}}\left(  \sum\limits_{j=1}^{n}\left\vert 1\right\vert
^{\delta^{p_{1}\cdots p_{m}}}\right)  ^{\frac{1}{\delta^{p_{1},\cdots,p_{m}}}%
}\right) \\
&  \leq n^{\frac{1}{\delta^{p_{1},\cdots,p_{m}}}}.
\end{align*}
On the other hand%
\[
\left(  \sum_{j_{1}=1}^{n}\left(  \sum_{j_{2}=1}^{n}\cdots\left(  \sum
_{j_{m}=1}^{n}B_{n}(e_{j_{1}},...,e_{j_{m}})^{q_{m}}\right)  ^{\frac{q_{m-1}%
}{q_{m}}}\cdots\right)  ^{\frac{q_{1}}{q_{2}}}\right)  ^{\frac{1}{q_{1}}%
}=n^{\frac{1}{q_{1}}},
\]
and, since $n$ is arbitrary,
\[
q_{1}\geq\delta^{p_{1},...,p_{m}}.
\]

Now let us prove the converse direction.

\bigskip We recall that for a bounded $m$-linear form $T:X_{p_{1}}\times
\cdots\times X_{p_{m}}\rightarrow\mathbb{K}$, we have%
\begin{align}
\left\Vert T\right\Vert  &  =\sup_{\left\Vert x^{\left(  i\right)
}\right\Vert _{X_{p_{i}}}\leq1;~1\leq i\leq m}\left\vert T\left(  x^{\left(
1\right)  },...,x^{\left(  m\right)  }\right)  \right\vert \label{nova1}\\
&  =\sup_{\left\Vert x^{\left(  i\right)  }\right\Vert _{X_{p_{i}}}%
\leq1;~1\leq i\leq m}\left\vert \sum_{i_{1},...,i_{m}=1}^{\infty}T\left(
e_{i_{1}},...,e_{i_{m}}\right)  x_{i_{1}}^{\left(  1\right)  }...x_{i_{m}%
}^{\left(  m\right)  }\right\vert \nonumber\\
&  =\sup_{\left\Vert x^{\left(  i\right)  }\right\Vert _{X_{p_{i}}}%
\leq1;~1\leq i\leq m-1}\left(  \sum_{i_{m}=1}^{\infty}\left\vert \sum
_{i_{1},...,i_{m-1}=1}^{\infty}T\left(  e_{i_{1}},...,e_{i_{m}}\right)
x_{i_{1}}^{\left(  1\right)  }...x_{i_{m-1}}^{\left(  m-1\right)  }\right\vert
^{\left(  p_{m}\right)  ^{\ast}}\right)  ^{\frac{1}{\left(  p_{m}\right)
^{\ast}}}.\nonumber
\end{align}
We denote by $X_{r}^{+}$ the set of sequences $\left(  x_{j}\right)  \in
X_{r}$, such that $x_{j}\geq0$ for all $j$. In the case $m=1$ the result is
immediate, it holds with constant $1$ and doesn't need the non-negative
assumption. Let us show the general case $m$, supposing that the result holds
for $m-1;$ so we suppose that if $p_{1},...,p_{m-1}\in(1,\infty]$ are such
that $\frac{1}{p_{1}}+\cdots+\frac{1}{p_{m-1}}<1$, then
\[
\left(  \sum_{j_{_{1}}=1}^{\infty}\left(  \sum_{j_{2}=1}^{\infty}\cdots\left(
\sum_{j_{_{m-1}}=1}^{\infty}A(e_{j_{1}},...,e_{j_{m-1}})^{\delta^{p_{m-1}}%
}\right)  ^{\frac{\delta^{p_{m-2},p_{m-1}}}{\delta^{p_{m-1}}}}\cdots\right)
^{\frac{\delta^{p_{1},...,p_{m-1}}}{\delta^{p_{2},...,p_{m-1}}}}\right)
^{\frac{1}{\delta^{p_{1},...,p_{m-1}}}}\leq\left\Vert A\right\Vert
\]
for all bounded non negative $\left(  m-1\right)  $-linear forms $A:X_{p_{1}%
}\times\cdots\times X_{p_{m-1}}\rightarrow\mathbb{K}.$

Suppose that $p_{1},...,p_{m}\in(1,\infty]$ are such that $\frac{1}{p_{1}%
}+\cdots+\frac{1}{p_{m}}<1.$ In this case%
\[
\frac{1}{p_{1}}+\cdots+\frac{1}{p_{m-1}}<1-\frac{1}{p_{m}}=\left(
\delta^{p_{m}}\right)  ^{-1}%
\]
and then for all $i$ $\in\left\{  1,..,m-1\right\}  $, we have $p_{i}%
\geq\delta^{p_{m}}$.

Let $D:X_{p_{1}}\times\cdots\times X_{p_{m}}\rightarrow\mathbb{K}$ be a
bounded non negative $m$-linear form. We define the bounded non negative
$\left(  m-1\right)  $-linear form $A:X_{r_{1}}\times\cdots\times X_{r_{m-1}%
}\rightarrow\mathbb{K}$ by%
\begin{equation}
A(e_{j_{1}},...,e_{j_{m-1}})=\sum_{j_{_{m}}=1}^{\infty}D(e_{j_{1}%
},...,e_{j_{m}})^{\delta^{p_{m}}}, \label{serie}%
\end{equation}
with $r_{i}=p_{i}/\delta^{p_{m}}$ for each $i\in\left\{  1,..,m-1\right\}  $.
Note that $A$ is well defined. In fact:

(i)%
\[
\sum_{j_{_{m}}=1}^{\infty}D(e_{j_{1}},...,e_{j_{m}})^{\delta^{p_{m}}}%
=\sum_{j_{_{m}}=1}^{\infty}\left\vert D(e_{j_{1}},...,e_{j_{m}})\right\vert
^{\left(  p_{m}\right)  ^{\ast}}<\infty,
\]
since $D(e_{j_{1}},...,e_{j_{m-1},}\cdot):X_{p_{m}}\rightarrow\mathbb{K}$ is a
linear form and $\left(  e_{j}\right)  $ is weakly $\left(  p_{m}\right)
^{\ast}$-summable.

(ii) For all positive integers $n,$ we have
\begin{align*}
&  \sup_{\left\Vert x^{\left(  i\right)  }\right\Vert _{X_{r_{i}}}\leq1;~1\leq
i\leq m-1}\sum_{j_{1},...,j_{m-1}=1}^{n}\left\vert A(e_{j_{1}},...,e_{j_{m-1}%
})x_{j_{1}}^{\left(  1\right)  }...x_{j_{m-1}}^{\left(  m-1\right)
}\right\vert \\
&  =\sup_{\left\Vert x^{\left(  i\right)  }\right\Vert _{X_{r_{i}}^{+}}%
\leq1;~1\leq i\leq m-1}\sum_{j_{1},...,j_{m-1}=1}^{n}A(e_{j_{1}}%
,...,e_{j_{m-1}})x_{j_{1}}^{\left(  1\right)  }...x_{j_{m-1}}^{\left(
m-1\right)  }\\
&  =\sup_{\left\Vert x^{\left(  i\right)  }\right\Vert _{X_{p_{i}}^{+}}%
\leq1;~1\leq i\leq m-1}\sum_{j_{1},...,j_{m-1}=1}^{n}A(e_{j_{1}}%
,...,e_{j_{m-1}})\left(  x_{j_{1}}^{\left(  1\right)  }...x_{j_{m-1}}^{\left(
m-1\right)  }\right)  ^{\delta^{p_{m}}}\\
&  =\sup_{\left\Vert x^{\left(  i\right)  }\right\Vert _{X_{p_{i}}^{+}}%
\leq1;~1\leq i\leq m-1}\sum_{j_{m}=1}^{\infty}\sum_{j_{1},...,j_{m-1}=1}%
^{n}D(e_{j_{1}},...,e_{j_{m}})^{\delta^{p_{m}}}\left(  x_{j_{1}}^{\left(
1\right)  }...x_{j_{m-1}}^{\left(  m-1\right)  }\right)  ^{\delta^{p_{m}}}\\
&  =\sup_{\left\Vert x^{\left(  i\right)  }\right\Vert _{X_{p_{i}}^{+}}%
\leq1;~1\leq i\leq m-1}\sum_{j_{m}=1}^{\infty}\left(  \sum_{j_{1}%
,...,j_{m-1}=1}^{n}D(e_{j_{1}},...,e_{j_{m}})x_{j_{1}}^{\left(  1\right)
}...x_{j_{m-1}}^{\left(  m-1\right)  }\right)  ^{\delta^{p_{m}}}%
\end{align*}
and by (\ref{nova1}) we conclude that%
\begin{equation}
\sum_{j_{1},...,j_{m-1}=1}^{n}\left\vert A(e_{j_{1}},...,e_{j_{m-1}})x_{j_{1}%
}^{\left(  1\right)  }...x_{j_{m-1}}^{\left(  m-1\right)  }\right\vert
\leq\left\Vert D\right\Vert ^{\delta^{p_{m}}} \label{2211}%
\end{equation}
for all $n$ and all $x^{(k)}\in B_{X_{r_{k}}},$ with $k=1,...,m-1$, and we
conclude that $A:X_{r_{1}}\times\cdots\times X_{r_{m-1}}\rightarrow\mathbb{K}$
is well defined.

Note that
\[
\frac{1}{r_{i}}+\cdots+\frac{1}{r_{m-1}}=\left(  \frac{1}{p_{i}}+\cdots
+\frac{1}{p_{m-1}}\right)  \delta^{p_{m}}=\left(  \frac{1}{p_{i}}+\cdots
+\frac{1}{p_{m-1}}\right)  \left(  1-\frac{1}{p_{m}}\right)  ^{-1}<1
\]
Hence, for each $i\in\left\{  1,..,m-1\right\}  ,$ a simple calculation shows
that%
\begin{equation}
\delta^{r_{i},...,r_{m-1}}=\frac{\delta^{p_{i},...,p_{m}}}{\delta^{p_{m}}}.
\label{nova2}%
\end{equation}
Therefore, by (\ref{serie}) and (\ref{nova2}) we have
\begin{align*}
&  \left(  \sum_{j_{_{1}}=1}^{\infty}\left(  \sum_{j_{_{2}}=1}^{\infty}%
\cdots\left(  \sum_{j_{_{m}}=1}^{\infty}D(e_{j_{1}},...,e_{j_{m}}%
)^{\delta^{p_{m}}}\right)  ^{\frac{\delta^{p_{m-1},p_{m}}}{\delta^{p_{m}}}%
}\cdots\right)  ^{\frac{\delta^{p_{1},...,p_{m}}}{\delta^{p_{2},...,p_{m}}}%
}\right)  ^{\frac{1}{\delta^{p_{1},...,p_{m}}}\times\delta^{p_{m}}}\\
&  =\left(  \sum_{j_{_{1}}=1}^{\infty}\left(  \sum_{j_{_{2}}=1}^{\infty}%
\cdots\left(  \sum_{j_{_{m-1}}=1}^{\infty}A(e_{j_{1}},...,e_{j_{m-1}}%
)^{\frac{\delta^{p_{m-1},p_{m}}}{\delta^{p_{m}}}}\right)  ^{\frac
{\delta^{p_{m-2},p_{m}}}{\delta^{p_{m-1},p_{m}}}}\cdots\right)  ^{\frac
{\delta^{p_{1},...,p_{m}}}{\delta^{p_{2},...,p_{m}}}}\right)  ^{\frac
{1}{\delta^{p_{1},...,p_{m}}}\times\delta^{p_{m}}}\\
&  =\left(  \sum_{j_{_{1}}=1}^{\infty}\left(  \sum_{j_{_{2}}=1}^{\infty}%
\cdots\left(  \sum_{j_{_{m-1}}=1}^{\infty}A(e_{j_{1}},...,e_{j_{m-1}}%
)^{\delta^{r_{m-1}}}\right)  ^{\frac{\delta^{r_{m-2,}r_{m-1}}}{\delta
^{r_{m-1}}}}\cdots\right)  ^{\frac{\delta^{r_{1},...,r_{m-1}}}{\delta
^{r_{2},...,r_{m}-1}}}\right)  ^{\frac{1}{\delta^{r_{1},...,r_{m}-1}}}.
\end{align*}
By the last equality and the Induction Hypothesis we conclude that
\begin{align*}
&  \left(  \sum_{j_{_{1}}=1}^{\infty}\left(  \sum_{j_{_{2}}=1}^{\infty}%
\cdots\left(  \sum_{j_{_{m}}=1}^{\infty}D(e_{j_{1}},...,e_{j_{m}}%
)^{\delta^{p_{m}}}\right)  ^{\frac{\delta^{p_{m-1},p_{m}}}{\delta^{p_{m}}}%
}\cdots\right)  ^{\frac{\delta^{p_{1},...,p_{m}}}{\delta^{p_{2},...,p_{m}}}%
}\right)  ^{\frac{1}{\delta^{p_{1},...,p_{m}}}\times\delta^{p_{m}}}\\
&  \leq\sup_{\left\Vert x^{\left(  i\right)  }\right\Vert _{X_{r_{i}}}%
\leq1;~1\leq i\leq m-1}\left\vert \sum_{j_{1},...,j_{m-1}=1}^{\infty
}A(e_{j_{1}},...,e_{j_{m-1}})x_{j_{1}}^{\left(  1\right)  }...x_{j_{m-1}%
}^{\left(  m-1\right)  }\right\vert \\
&  \leq\left\Vert D\right\Vert ^{\delta^{p_{m}}},
\end{align*}
where in the last inequality we have used (\ref{2211}).

\bigskip\textbf{Second Case. }$\frac{1}{p_{1}}+\cdots+\frac{1}{p_{m}}\geq1.$

\bigskip We begin by proving the direct implication.

Consider%
\[
A(x^{(1)},...,x^{(m)})=\sum\limits_{j=1}^{\infty}x_{j}^{(1)}\cdots x_{j}%
^{(m)},
\]
and we conclude that%
\[
q_{1}=\infty=\delta^{p_{1},...,p_{m}}.
\]
If
\[
\frac{1}{p_{i}}+\cdots+\frac{1}{p_{m}}\geq1
\]
for all $i$, the proof is immediate. Otherwise, at some stage $(i=2,3,...)$ we
begin to have a strict inequality
\[
\frac{1}{p_{i}}+\cdots+\frac{1}{p_{m}}<1.
\]
Denote by $i_{0}$ this index. Then
\[
\frac{1}{p_{i_{0}-1}}+\cdots+\frac{1}{p_{m}}\geq1
\]
and%
\[
\frac{1}{p_{i_{0}}}+\cdots+\frac{1}{p_{m}}<1.
\]
If $i_{0}=3$, we consider
\[
A(x^{(1)},...,x^{(m)})=x_{1}^{(1)}\sum\limits_{j=1}^{\infty}x_{j}^{(2)}\cdots
x_{j}^{(m)}%
\]
and we conclude that%
\[
q_{1}=q_{2}=\infty=\delta^{p_{2},...,p_{m}}.
\]
Similarly, if $i_{0}>3$, we have%
\[
q_{1}=q_{2}=...=q_{i_{0}-1}=\infty=\delta^{p_{i_{0}-1},...,p_{m}}%
\]
and%
\[
\sup_{j_{1},..,j_{i_{0}-1}}\left(  \sum_{j_{i_{0}}=1}^{\infty}\left(
\cdots\left(  \sum_{j_{m}=1}^{\infty}A(e_{j_{1}},...,e_{j_{m}})^{q_{m}%
}\right)  ^{\frac{q_{m-1}}{q_{m}}}\cdots\right)  ^{\frac{q_{i_{0}+1}}%
{q_{i_{0}}}}\right)  ^{\frac{1}{q_{i_{0}}}}\leq\left\Vert A\right\Vert .
\]

Now it is simple to imitate the arguments of the previous case to complete the proof.

\bigskip Now we prove the reverse implication. Recall that we are in the case
\[
\frac{1}{p_{1}}+\cdots+\frac{1}{p_{m}}\geq1.
\]
If
\[
\frac{1}{p_{i}}+\cdots+\frac{1}{p_{m}}\geq1
\]
for all $i$, the proof is immediate. Otherwise, at some stage $(i=2,3,...)$ we
begin to have a strict inequality
\[
\frac{1}{p_{i}}+\cdots+\frac{1}{p_{m}}<1.
\]
Denote by $i_{0}$ this index. Then
\[
\frac{1}{p_{i_{0}-1}}+\cdots+\frac{1}{p_{m}}\geq1
\]
and%
\[
\frac{1}{p_{i_{0}}}+\cdots+\frac{1}{p_{m}}<1.
\]
We need to prove that%
\[
\sup_{j_{1},..,j_{i_{0}-1}}\left(  \sum_{j_{i_{0}}=1}^{\infty}\left(
\cdots\left(  \sum_{j_{m}=1}^{\infty}A(e_{j_{1}},...,e_{j_{m}})^{q_{m}%
}\right)  ^{\frac{q_{m-1}}{q_{m}}}\cdots\right)  ^{\frac{q_{i_{0}}}%
{q_{i_{0}+1}}}\right)  ^{\frac{1}{q_{i_{0}}}}\leq\left\Vert A\right\Vert
\]
for
\[
q_{i_{0}}\geq\delta^{p_{i_{0}},...,p_{m}},...,q_{m}\geq\delta^{p_{m}}.
\]
By the first case, we know that for any fixed vectors $e_{j_{1}}%
,...,e_{j_{i_{0}-1}}$, we have
\[
\left(  \sum_{j_{i_{0}}=1}^{\infty}\left(  \cdots\left(  \sum_{j_{m}%
=1}^{\infty}A(e_{j_{1}},...,e_{j_{m}})^{\delta^{p_{m}}}\right)  ^{\frac
{\delta^{p_{m-1},p_{m}}}{\delta^{p_{m}}}}\cdots\right)  ^{\frac{\delta
^{p_{i_{0}},...,p_{m}}}{\delta^{p_{i_{0}+1},...,p_{m}}}}\right)  ^{\frac
{1}{\delta^{p_{i_{0}},...,p_{m}}}}\leq\left\Vert A\right\Vert
\]
for all bounded non negative $m$-linear forms $A\colon X_{p_{1}}\times
\cdots\times X_{p_{m}}\rightarrow\mathbb{K}$. Then,
\[
\sup_{j_{1},..,j_{i_{0}-1}}\left(  \sum_{j_{i_{0}}=1}^{\infty}\left(
\cdots\left(  \sum_{j_{m}=1}^{\infty}A(e_{j_{1}},...,e_{j_{m}})^{\delta
^{p_{m}}}\right)  ^{\frac{\delta^{p_{m-1},p_{m}}}{\delta^{p_{m}}}}%
\cdots\right)  ^{\frac{\delta^{p_{i_{0}},...,p_{m}}}{\delta^{p_{i_{0}%
+1},...,p_{m}}}}\right)  ^{\frac{1}{\delta^{p_{i_{0}},...,p_{m}}}}%
\leq\left\Vert A\right\Vert
\]
for all bounded non negative $m$-linear forms $A\colon X_{p_{1}}\times
\cdots\times X_{p_{m}}\rightarrow\mathbb{K}$.

\end{document}